\documentclass[12pt]{article} 
\usepackage[cp1251]{inputenc}
\usepackage{graphicx}
\usepackage{cite}
\usepackage{textcase}
\usepackage{amsfonts}
\usepackage{amsmath}
\usepackage[russian]{babel}

\newtheorem{theorem}{Теорема}
\newtheorem{remark}{Замечание}

\begin{document} 

\begin{center}
\textbf{Нелинейные законы управления, построенные на базе линейных с использованием нечетных функций}

\end{center}

\begin{center}
 Игорь Борисович~Фуртат
\end{center}

\begin{center}
Институт проблем машиноведения РАН, 2021 г.
\end{center}

\begin{abstract}
В статье исследуются нелинейные законы управления, полученные из линейного путем двух типов замен с использованием нечетных функций. 
Первая замена состоит в пропускании каждой компоненты вектора состояния через нелинейную функцию,  вторая замена -- в пропускании всего линейного закона управления через нелинейную функцию. 
Для исследования таких систем предлагается нелинейные функции представить в виде линейных с нелинейным угловым коэффициентом.
Такое представление позволит использовать аппарат линейных матричных неравенств (ЛМН) для исследования устойчивости замкнутых систем. 
Найдены области устойчивости и области для начальных условий, полученные в результате многошаговой процедуры поиска решений с использованием ЛМН.
Показано, что использование предложенных нелинейных законов управления позволяет уменьшить установившуюся ошибку по сравнению с линейным.

\end{abstract}

\textit{Ключевые слова:} динамическая система, нелинейный закон упрвления, линейные матричные неравенства.

\section{Введение.}

В теорие автоматического управления существует практика использования различных нелинейных функций, замещающих или дополняющих существующие законы управления. 
Например, знаковые законы управления аппроксимируют сигмоидальными функциями  \cite{Khalil09,Krasnova15} для реализации на практике, а также возможности избавиться от чаттеринга в окрестности плоскости скольжения и исследования замкнутой системы без разрывных нелинейностей. 
В \cite{Krasnova16,Antipov21} ранее разработанный закон управления \cite{Krasnova06} пропускается через функцию гиперболического тангенса для уменьшения установившейся ошибки слежения по сравнению с \cite{Krasnova06}. 
Также в \cite{Khalil09,Krasnova15,Krasnova16,Antipov21} отмечается, что помимо повышения качества регулирования в установившемся режиме, использование сигмоидальных функций позволяет сформировать заранее ограниченный закон управления.

Однако ограниченный закон управления в условиях возмущений и неустойчивого объекта может повлечь  рассмотрение ограничений на фазовое пространство для сохранения устойчивости замкнутой системы. 
В условиях же ограниченного фазового пространства можно рассматривать любые непрерывные функции, не только сигмоидальные. 

В \cite{Khalil09} показано, что аддитивная добавка, например, функции $x^3$ к существующему закону управления повышает качество демпфирования внешних возмущений. 
В адаптивном управлении \cite{Andrievkiy00} суммирование сигмоидальных слагаемых к существующему адаптивному алгоритму повышает качество регулирования в установившемся режиме.
В \cite{Nekhoroshikh20a,Nekhoroshikh20b} замена каждой компоненты вектора состояния $x_i$ в линейном законе управления на соответствующие функции $x_i^{\gamma} \textup{sgn}\{x_i\}$, $\gamma>0$ позволяет увеличить скорость сходимости нормы вектора состояния замкнутой системы по сравнению с использованием линейного алгоритма.

Однако в \cite{Khalil09,Krasnova15,Krasnova16,Antipov21,Nekhoroshikh20a,Nekhoroshikh20b} использование нелинейных функций зачастую требует перерасчета настроечных коэффициентов для нового закона управления, что может усложнить процесс синтеза алгоритма. 
Проблема же введения нелинейных функций, не меняющих настроечных коэффициентов в существующих законах управления, но при этом улучшающих качество управления в установившемся режиме, в работах \cite{Khalil09,Krasnova15,Krasnova16,Antipov21,Nekhoroshikh20a,Nekhoroshikh20b} не рассматривалась.
Поэтому в данной статье на базе линейного закона управления рассмотрим построение новых нелинейных законов, не меняя настроечных параметров исходного линейного алгоритма. 
Также предложим новый метод исследования устойчивости таких систем, основанный на представлении нелинейной функции в линейной форме с нелинейным угловым коэффициентом. 
Такая форма позволит использовать аппарат линейных матричных неравенств (ЛМН) и упростить процедуру поиска множества устойчивости и множества начальных условий. 
Предложенная в статье многошаговая процедура поиска решений на базе ЛМН позволит расширить искомые области устойчивости и начальных условий.

Статья организована следующим образом. 
В разделе \ref{Sec2} задается линейный закон управления, стабилизирующий линейную систему. 
Не меняя настроечных коэффициентов в линейном законе, каждая компонента вектора состояния пропускается через нечетную функцию, либо весь линейный закон управления пропускается через нечетную функцию. 
В разделе \ref{Sec3} описываются предложенные методы решения и формулируются условия на свойства нелинейных функций, множества устойчивости и множества начальных условий при выполнении которых нелинейные законы управления будут стабилизировать объект управления.
В разделе \ref{Sec4} показано, что предложенные нелинейные алгоритмы управления гарантируют выше точность в установившемся режиме, чем линейный закон управления.

\section{Постановка задачи.}
\label{Sec2}

Рассмотрим линейный объект
\begin{equation}
\label{eq_system}
\begin{array}{l}
\dot{x}=Ax+Bu+Df,
\end{array}
\end{equation}
где $x = col\{x_1,...,x_n\}$ -- вектор состояния, 
$u \in \mathbb R$ -- сигнал управления, 
$f \in \mathbb R^l$ -- неизвестное возмущение такое, что $|f| \leq \bar{f}$, 
матрицы $A$, $B$ и $D$ известны и имеют соответствующие размерности,
пара $(A,B)$ управляема. 

Известно \cite{Andrievkiy00}, что существует матрица $K=[k_1,...,k_n]$ такая, что закон управления 
\begin{equation}
\label{eq_known_CL}
\begin{array}{l}
u = Kx = k_1 x_1 + ... + k_n x_n
\end{array}
\end{equation}
гарантирует выполнение предельного неравенства
\begin{equation}
\label{eq_goal}
\begin{array}{l}
\overline{\lim\limits_{t \to \infty}}|x(t)| \leq \delta.
\end{array}
\end{equation}
Здесь $\delta>0$ -- точность регулирования в установившемся режиме.

В статье введем в рассмотрение нечетную функцию $\varphi(\cdot)$. 
С учетом данной функции и \eqref{eq_known_CL}, построим два новых закона управления:
\begin{equation}
\label{eq_sigma_CL}
\begin{array}{l}
u = k_1 \varphi(x_1) + ... + k_n \varphi(x_n),
\end{array}
\end{equation}
\begin{equation}
\label{eq_sigma_CL_all}
\begin{array}{l}
u = \varphi \left( k_1 x_1 + ... + k_n x_n \right).
\end{array}
\end{equation}

Требуется определить условия на функцию $\varphi(\cdot)$, при которых при заданной матрице $K$ в \eqref{eq_known_CL} законы управления \eqref{eq_sigma_CL} и \eqref{eq_sigma_CL_all} обеспечат выполнение предельного неравенства \eqref{eq_goal} с величиной $\delta$ меньше, чем линейный закон управления \eqref{eq_known_CL}.

\begin{remark}
\label{Rem1}
В качетсве $\varphi(\cdot)$ могут быть рассмотрены стандартные линейные, степенные и сигмоидальные функции, функция насыщения (см. рис. \ref{Fig_sigma_Gen}) и т.п. 
Другие специфические функции будут приведены далее.

\begin{figure}[h!]
\center{\includegraphics[width=0.6\linewidth]{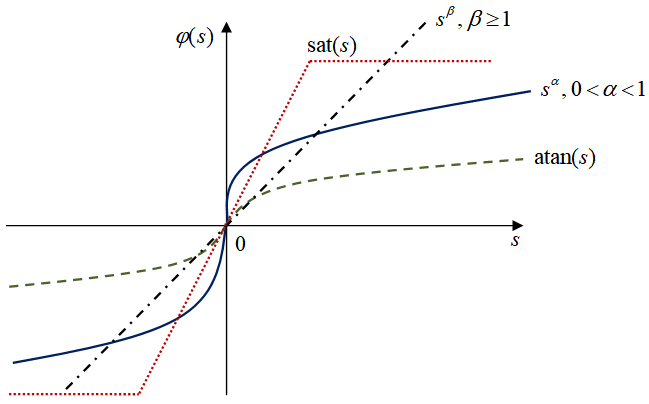}}
\caption{Примеры функций $\varphi(s)$.}
\label{Fig_sigma_Gen}
\end{figure}

\end{remark}

\section{Метод решения.}
\label{Sec3}

\begin{figure}[h!]
\center{\includegraphics[width=0.6\linewidth]{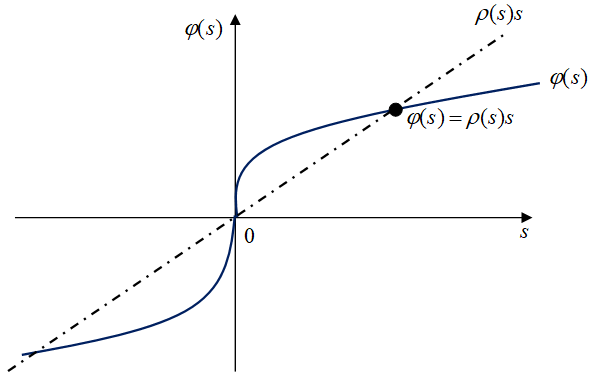}}
\caption{Иллюстрация замены функции $\varphi(s)$ на функцию линейного вида с переменным угловым коэффициентом $\rho(s)$.}
\label{Rho}
\end{figure}

Для исследования законов управления \eqref{eq_sigma_CL} и \eqref{eq_sigma_CL_all}, рассмотрим представление нелинейной функции $\varphi(s)$ в виде линейной $\varphi(s)=\rho(s) s$ с нелинейным угловым коэффициентом $\rho(s)$,  $s \in \mathbb R$ (см. рис. \ref{Rho}). 
Такой подход позволит исследовать нелинейные законы управления \eqref{eq_sigma_CL}, \eqref{eq_sigma_CL_all} наряду с линейным \eqref{eq_known_CL}, а также привлечь аппарат линейных матричных неравенств для анализа устойчивости замкнутой системы.

\subsection{Анализ устойчивости замкнутой системы с законом управления \eqref{eq_sigma_CL}.}
\label{SubSec2.1}

Обозначим $\Psi(x)=\textup{diag}\{\rho(x_1),...,\rho(x_n)\}$. 
Подставив \eqref{eq_sigma_CL} в \eqref{eq_system}, запишем уравнение замкнутой системы
\begin{equation}
\label{eq_closed_loop}
\begin{array}{l}
\dot{x} = (A+BK\Psi(x))x + Df.
\end{array}
\end{equation}


\begin{theorem}
\label{Th1}
Пусть для заданных $K$ в \eqref{eq_known_CL} и $\tau_i>0$ существуют $0 \leq \underline{\rho}=\rho_0<\rho_1<...<\rho_v=\bar{\rho}$, $\chi_i>0$, $\gamma_i>0$ и $P_i=P_i^{\textrm T}>0$ такие, что выполнены следующие матричные неравенства
\begin{equation}
\label{eq_LMI_P}
\begin{array}{l}
\begin{bmatrix}
(A+BK \Psi_i)^{\textrm T} P_i + P_i (A+BK \Psi_i)+\tau P_i & P_i D \\
* & -\chi_i I
\end{bmatrix}
< 0,
~~~
P_i > \gamma_i I
\end{array}
\end{equation} 
в вершинах $\Psi_i=\textup{diag}\{\{\rho_{i-1},\rho_{i}\},...,\{\rho_{i-1},\rho_{i}\}\}$, $i=1,...,v$.
Тогда область устойчивости $\mathcal X$, 
множество допустимых начальных условий $\mathcal X_0$ и верхняя оценка времени $T$ вхождения в инвариантную область $|x(t)|<\varepsilon$ определятся в виде
\begin{equation}
\label{eq_lim_XXX}
\begin{array}{l}
\mathcal X = \left \{ x \in \mathbb R^n: \underline{x} \leq |x_i| \leq \bar{x},~ i=1,...,n  \right \},
\end{array}
\end{equation}

\begin{equation}
\label{eq_Lyap_XXX0}
\begin{array}{l}
\mathcal X_0 = 
\left\{ 
x \in \mathbb R^n: 
|x| \leq  \bar{x}_0:=
\sqrt{
\frac{\underline{\tau} \underline{\gamma} \bar{x}^2 - 2 \bar{\chi} \bar{f}}{n^2 \underline{\tau} \|\bar{P}\|}
}
\right\},
\end{array}
\end{equation}

\begin{equation}
\label{eq_T_ult_bound}
\begin{array}{l}
T = 
\frac{1}{\underline{\tau}}
\ln
\frac
{\underline{\tau} |x(0)|^2 \|\bar{P}\| +  \bar{\chi} \bar{f}}
{\underline{\tau} \underline{\gamma} \varepsilon^2 - \bar{\chi} \bar{f}},
\end{array}
\end{equation}
$\underline{x}>0$ и $\bar{x}>0$ находятся как 
$s \in [-\bar{x}; -\underline{x}] \cup [\underline{x}; \bar{x}]$ из решения неравенств 
$\varphi(s)-\underline{\rho} s \geq 0$ и 
$\varphi(s)-\bar{\rho} s \leq 0$, 
$\underline{\tau}=\min\limits_{i=1,...,v}\{\tau_i\}$, 
$\underline{\gamma}=\min\limits_{i=1,...,v}\{\gamma_i\}$, 
$\bar{\chi}=\max\limits_{i=1,...,v}\{\chi_i\}$ и 
$\|\bar{P}\|=\max\limits_{i=1,...,v}\{\|P_i\| \}$.

\end{theorem}



\textbf{Доказательство.} 
Для анализа устойчивости замкнутой системы выберем функции Ляпунова
\begin{equation}
\label{eq_Lyap_fun}
\begin{array}{l}
V_i=x^{\textrm T} P_i x, ~~~ i=1,...,v
\end{array}
\end{equation}
и потребуем выполнение условий
\begin{equation}
\label{eq_Lyap_ineq}
\begin{array}{l}
\dot{V}_i + \tau_i V_i - \chi_i f^{\textrm T} f < 0, ~~~ i=1,...,v.
\end{array}
\end{equation}

Подставляя \eqref{eq_closed_loop} и \eqref{eq_Lyap_fun} в \eqref{eq_Lyap_ineq}, получим
\begin{equation}
\label{eq_Lyap_ineq_2}
\begin{array}{l}
x^{\textrm T}
\left[ (A+BK\Psi_i(x))^{\textrm T}P_i + P_i(A+BK\Psi_i(x)) \right]x +
\\
+ 2x^{\textrm T}P_iDf
+ \tau_i x^{\textrm T} P_i x 
- \chi_i f^{\textrm T} f < 0, ~~~ i=1,...,v.
\end{array}
\end{equation}

Введя новый вектор $z=col\{x, f\}$, перепишем \eqref{eq_Lyap_ineq_2} как
\begin{equation}
\label{eq_Lyap_ineq_3}
\begin{array}{l}
z^{\textrm T}
\begin{bmatrix}
(A+BK\Psi_i(x))^{\textrm T}P_i + P_i(A+BK\Psi_i(x))+\tau P & P_iD \\
* & -\chi_i I
\end{bmatrix}
z< 0, ~~~ i=1,...,v.
\end{array}
\end{equation}

Неравенство \eqref{eq_Lyap_ineq_3} будет выполнено, если будет выполнено условие 
\begin{equation}
\label{eq_LMI_proof_1}
\begin{array}{l}
\begin{bmatrix}
(A+BK\Psi_i(x))^{\textrm T}P_i + P_i(A+BK\Psi_i(x))+\tau P & P_iD \\
* & -\chi_i I
\end{bmatrix}
< 0, ~~~ i=1,...,v.
\end{array}
\end{equation}

Элементы диагональной матрицы $\Psi_i(x)$ принадлежат отрезку $[\rho_{i-1}; \rho_{i}]$. 
Значит, ЛМН \eqref{eq_LMI_proof_1} содержит политопную неопределенность. 
Согласно лемме \cite{11}, такое неравенство будет выполнено, если оно будет выполнено в вершинах политопа $[\rho_{i-1}; \rho_{i}]$. 
Значит, замкнутая система \eqref{eq_closed_loop} устойчива.


Пусть для любых $\underline{x}>0$ и $\bar{x}>0$ выполнены неравенства $\varphi(s)-\underline{\rho} s \geq 0$ и $\varphi(s)-\bar{\rho} s \leq 0$. 
Тогда $\underline{x}$ и $\bar{x}$ являются нижней и верхней оценками множества устойчивости, записанного в виде \eqref{eq_lim_XXX}.





Запишем решение неравенства \eqref{eq_Lyap_ineq} в виде
\begin{equation}
\label{eq_Lyap_ineq_4}
\begin{array}{l}
\underline{\gamma} |x|^2 < x^{\textrm T} P_i x < \frac{\bar{\chi} \bar{f}}{\underline{\tau}} 
+\left( \|\bar{P}\| \bar{x}_0^2 + \frac{\bar{\chi} \bar{f}}{\underline{\tau}} \right) e^{-\underline{\tau} t}, 
~~~
i=1,...,v.
\end{array}
\end{equation}
Учитывая \eqref{eq_lim_XXX} и \eqref{eq_Lyap_ineq_4} и приняв
 $\underline{\gamma} n^2 \bar{x}^2 = \frac{2\bar{\chi} \bar{f}}{\underline{\tau}}
+\|\bar{P}\| \bar{x}_0^2$, получим результат \eqref{eq_Lyap_XXX0}.

Верхняя оценка времени $T$ при котором $|x(t)|$ попадут в заданную область $|x(t)| < \varepsilon$ и больше ее не покинут, определяется в виде \eqref{eq_T_ult_bound} из \eqref{eq_Lyap_ineq_4}.

\begin{remark}
\label{Rem3}
В формулировке теоремы разрешимость \eqref{eq_LMI_P} ищется во внутренних точках $\underline{\rho}=\rho_0<\rho_1<...<\rho_v=\bar{\rho}$ интервала, а не сразу на границах $\underline{\rho}$ и $\bar{\rho}$ согласно лемме \cite{Fridman14}.  
Это связано с тем, что применив лемму \cite{Fridman14} для поиска $\underline{\rho}$ и $\bar{\rho}$ при которых ЛМН \eqref{eq_LMI_P} будет разрешимо, можно получить более грубые решения, чем при использовании следующей многошаговой процедуры. 
Согласно данной процедуре предлагается сначала определить интервал $(\underline{\rho},\rho_1)$, в вершинах которого ЛМН \eqref{eq_LMI_P} имеет решение, воспользовавшись леммой \cite{Fridman14}. 
Затем, ищется следующий интервал $(\rho_1,\rho_2)$, в вершинах которого ЛМН \eqref{eq_LMI_P} имеет решение. 
Данная процедура проделывается до тех пор, пока не находится предельное значение $\rho_k=\bar{\rho}$ такое, что в вершинах $(\rho_{k-1},\rho_{k})$ ЛМН \eqref{eq_LMI_P} еще сохраняется решение. 
В результате применения многошаговой процедуры длина интервала $(\underline{\rho};\bar{\rho})$ может быть больше, чем при использовании только леммы \cite{Fridman14}.
\end{remark}

\begin{remark}
Для применения многошаговой процедуры (см. замечание \ref{Rem3}) можно воспользоваться следующими рекомендациями для поиска $\bar{\rho}$.
Если $\frac{d \varphi(s)}{d s} \Big|_{s=0}=\rho(0) < \infty$, то $\bar{\rho} \leq \rho(0)$ (см. рис. \ref{Solut_derivat} слева). 
Если $\frac{d \varphi(s)}{d s} \Big|_{s=0}=\infty$, то $\bar{\rho}$ увеличивается до тех пор, пока ЛМН \eqref{eq_LMI_P} еще будет иметь решение на интервале $(\rho_{v-1};\rho_{v})$  (см. рис. \ref{Solut_derivat} справа). В некоторых случая возможно $\bar{\rho}=\infty$, что соотвествует обратной связи с большим коэффициентом усиления.

\begin{figure}[h]
\begin{minipage}[h]{0.49\linewidth}
\center{\includegraphics[width=0.95\linewidth]{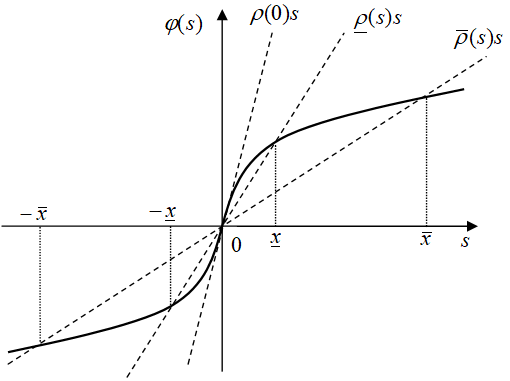}}
\end{minipage}
\hfill
\begin{minipage}[h]{0.49\linewidth}
\center{\includegraphics[width=0.95\linewidth]{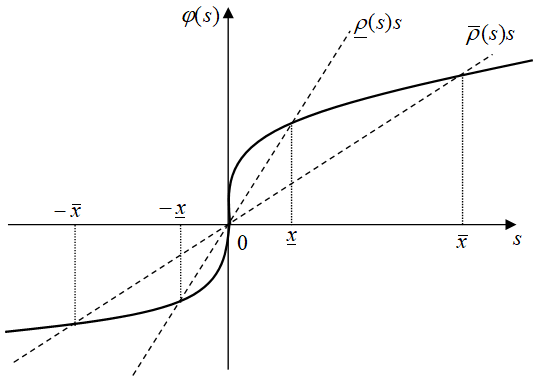}}
\end{minipage}
\caption{Поиск $\bar{\rho}$ при $\rho(0)<\infty$ (слева) и при $\rho(0)=\infty$ (справа).}
\label{Solut_derivat}
\end{figure}

\end{remark}

\begin{remark}
\label{Rem4}
Результаты секции \ref{Sec2} можно обобщить на закон управления вида 
$u = k_1 \varphi_1(x_1) + ... + k_n \varphi_n(x_n)$.
\end{remark}

\subsection{Анализ устойчивости замкнутой системы с законом управления \eqref{eq_sigma_CL_all}.}

Переписав $\varphi(KX)$ в виде $\varphi(Kx)=\rho(Kx) Kx$ и подставив \eqref{eq_sigma_CL_all} в \eqref{eq_system}, получим 
\begin{equation}
\label{eq_closed_loop_all}
\begin{array}{l}
\dot{x} = (A+BK\rho(Kx))x + Df.
\end{array}
\end{equation}


\begin{theorem}
\label{Th2}
Пусть для заданных $K$ и $\tau_i>0$ существуют $0 \leq \underline{\rho}=\rho_0<\rho_1<...<\rho_v=\bar{\rho}$, $\chi_i>0$, $\gamma_i>0$ и $P_i=P_i^{\textrm T}>0$ такие, что выполнены матричные неравенства 
\eqref{eq_LMI_P}
в вершинах $\Psi_i=\{\rho_{i-1},\rho_{i}\}$, $i=1,...,v$.
Тогда область устойчивости $\mathcal X$, 
множество допустимых начальных условий $\mathcal X_0$ и верхняя оценка времени переходного процесса $T$ в область $|x(t)|<\varepsilon$ определяются в виде
\eqref{eq_lim_XXX},
\eqref{eq_Lyap_XXX0},
\eqref{eq_T_ult_bound},
где $\underline{x}$ и $\bar{x}$ находятся как $s \in [\underline{x}; \bar{x}]$: 
$\varphi(s \sum_{i=1}^{n} k_i)-\underline{\rho} s \sum_{i=1}^{n} k_i \geq 0$ и $\varphi(s \sum_{i=1}^{n} k_i)-\bar{\rho} s \sum_{i=1}^{n} k_i \leq 0$.

\end{theorem}



\textbf{Доказательство.} 
Ради простоты положим $x_1=...=x_n=s$.  
Пусть для любого $s \in [\underline{x}; \bar{x}]$ выполнены неравенства $\varphi(s \sum_{i=1}^{n} k_i)-\underline{\rho} s \sum_{i=1}^{n} k_i \geq 0$ и $\varphi(s \sum_{i=1}^{n} k_i)-\bar{\rho} s \sum_{i=1}^{n} k_i \leq 0$. 
Тогда $\underline{x}$ и $\bar{x}$ являются нижней и верхней оценками множества устойчивости, записанного в виде \eqref{eq_lim_XXX}.

Поскольку структура \eqref{eq_closed_loop_all} подобна структуре \eqref{eq_closed_loop}, то дальнейшее доказательство теоремы \ref{Th2} аналогично доказательству теоремы \ref{Th1}.

\section{Анализ значения $\delta$ в \eqref{eq_goal} для линейного и нелинейных регуляторов и примеры функций $\varphi(s)$.}
\label{Sec4}

Покажем, что алгоритмы \eqref{eq_sigma_CL} и \eqref{eq_sigma_CL_all} обеспечивают лучше предельную ошибку, по сравнению с алгоритмом \eqref{eq_known_CL}. 
Для этого необходимо рассмотреть решение любой оптимизационной задачи. 
Например, рассмотрим оптимизационную задачу из \cite{Boyd94}, которая в нашем случае будет сформулирована в следующем виде.
Для заданных $K$ и $\tau$ требуется найти $P$ и $\gamma$, удовлетворяющие следующей оптимизационной задаче:
\begin{equation}
\label{eq_optim_task}
\begin{array}{l}
\textup{maximize}~\gamma
\\
\textup{subject to}~
\begin{bmatrix}
(A+BK \Theta)^{\textrm T}P + P(A+BK \Theta)+\tau P & PD \\
* & -\tau I
\end{bmatrix}
< 0
~\textup{and}
~P-\gamma I>0.
\end{array}
\end{equation}
Тогда 
\begin{equation}
\label{eq_ult_bound}
\begin{array}{l}
\overline{\lim}_{t \to \infty}|x(t)| \leq \delta := \sqrt{\frac{\bar{f}}{\gamma}}. 
\end{array}
\end{equation}

\begin{figure}[h]
\center{\includegraphics[width=0.5\linewidth]{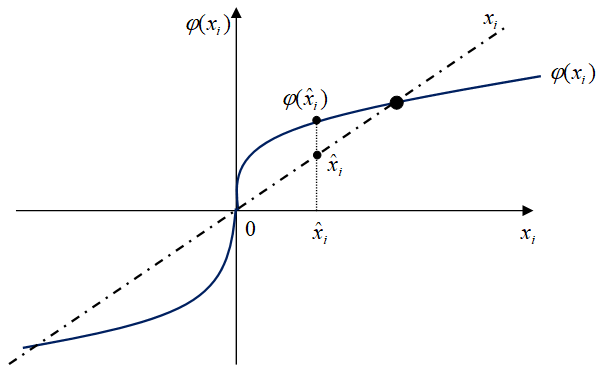}}
\caption{Сравнение линейного и нелинейного регуляторов.}
\label{Better}
\end{figure}

Для нелинейного регулятора \eqref{eq_known_CL} оптимизационная задача решается при $\Theta=I$, для нелинейного при $\Theta=\bar{\rho}I$ в установившемся режиме.
Поскольку для нелинейного регулятора $\|\Theta\|>1$ (см. рис. \ref{Better}, где для каждого фиксированного значения $\hat{x}_i$: $|\hat{x}_i|<|\varphi(\hat{x}_i)|$), значит, $\gamma$ для \eqref{eq_known_CL} будет меньше, чем для  \eqref{eq_sigma_CL} и \eqref{eq_sigma_CL_all} .
Другими словами, при  $\|\Theta\|>1$ происходит увеличение коэффициента в обратной связи, что делает матрицу $A+BK \Theta$ <<более>> устойчивой. 
Для <<более>> устойчивой матрицы требуется <<больше>> матрица $P$ (а, значит, и больше $\gamma$), что соотвествует <<меньшему>> притягивающему эллипсоиду в установившемся режиме.

Приведем некоторые примеры функций, для которых справедливы полученные результаты:
\begin{enumerate}

\item [(1)] $\varphi(s)=\mu \textup{sat}(\sigma s)$, $\varphi(s)=\mu \arctg(\sigma s)$, $\varphi(s)=\mu \frac{1-e^{-0.5\sigma s}}{1+e^{-0.5\sigma s}}$ и т.п. 
Достоинство таких функций состоит в том, что они гарантируют ограниченность закона управления. 
Однако для них всегда существует $\bar{x} < \infty$.

\item [(2)] $\varphi(s)=s^{\lambda}$, $0<\lambda<1$, $\varphi(s)=s^{\psi(s)}$ с четной, ограниченной, положительной функцией $\psi(s)$ такой, что $\psi(0) < 1$ и $\psi(\pm \infty) < \infty$ (например, $\psi(s)=\mu \frac{s^2+\mu^{-2}}{s^2+1}$, $\mu>0$), $\varphi(s)=s^{\lambda}+s^{1/\lambda}$ и т.п. Достоинство таких функций состоит в том, что они гарантируют ускоренную сходимость. Однако для них может существовать $\underline{x} < \infty$ и всегда существует $\bar{x} < \infty$.

\item [(3)] $\varphi(s)=\mu \textup{sat}(\sigma s) + \vartheta s$, $\varphi(s)=\mu \arctg(\sigma s) + \vartheta s$, $\varphi(s)=\mu \frac{1-e^{-0.5\sigma s}}{1+e^{-0.5\sigma s}} + \vartheta s$ и т.п. По сравнению с (1) данные функции неограниченные. 
Однако за счет соответствующего выбора $\mu$, $\sigma$ и $\vartheta$ можно обеспечить $\underline{x}=0$ и $\bar{x}=\infty$.

\end{enumerate}


\end{document}